\documentclass[a4paper]{amsart}
\usepackage[T1]{fontenc}
\usepackage[english]{babel}
\usepackage[cp1252]{inputenc}
\usepackage{amsthm}
\usepackage{amsmath}
\usepackage{amsfonts}
\usepackage{stmaryrd}
\usepackage{amssymb}
\usepackage{hyperref}
\usepackage{epsfig,epsf}
\usepackage{indentfirst}
\usepackage{amssymb}
\usepackage{eufrak}
\usepackage{mathrsfs}
\usepackage{xypic}

\newcommand{\Verte}{\textrm Vert}
\newcommand{\di}{\textrm{d}}
\newcommand{\R}{\mathbb R}
\newcommand{\C}{\mathbb C}
\newcommand{\Log}{\mathrm{Log}\,}
\newcommand{\A}{\mathscr A}
\newcommand{\Si}{\mathscr S}
\newcommand{\de}{\textrm{deg}}

\newtheorem{remark}{Remark}[section]

\newtheorem{theorem}{Theorem}[section]

\newtheorem{proposition}{Proposition}[section]

\newtheorem{lemma}{Lemma}[section]

\begin{document}
\title{On the volume of complex amoebas}
\author{Farid Madani and Mounir Nisse}

\date{}
\address{NWF I-Mathematik, Universit\"at Regensburg, 93040 Regensburg, Germany.}
\email{\href{mailto:Farid.Madani@mathematik.uni-regensburg.de}{Farid.Madani@mathematik.uni-regensburg.de}}
\thanks{The first author is supported by the Alexander von Humboldt Foundation.}

\address{Department of Mathematics, Texas A\&M University, College Station, TX 77843-3368, USA.}
\email{\href{mailto:nisse@math.tamu.edu}{nisse@math.tamu.edu}}
\thanks{The research of the second author is partially supported by NSF MCS grant DMS-0915245.}
\maketitle
\vspace{-0.5cm}
\begin{center}\em{To the memory of Mikael Passare}\end{center}
\begin{abstract} The paper deals with amoebas of  $k$-dimensional algebraic varieties  in the   complex algebraic torus of dimension $n\geq 2k$. First, we show that the area of  complex algebraic curve amoebas is finite. Moreover, we give an estimate of this area in the rational curve case in terms of the  degree  of the  rational parametrization
coordinates. We also show that the volume of  the amoeba of a  $k$-dimensional algebraic variety  in $(\mathbb{C}^*)^{n}$,  with $n\geq 2k$,  is finite.
\end{abstract}


\section{Introduction}

Amoebas have proven to be a very useful tool in several areas of mathematics, and  they have many applications in real algebraic geometry, complex analysis, mirror symmetry, algebraic statistics and in several other areas (see \cite{M1-02}, \cite{M2-04}, \cite{M3-00}, \cite{FPT-00}, \cite{NS-11}, \cite{PR-04}, \cite{PS-04} and \cite{R-01}). They degenerate to  piecewise-linear objects called {\em tropical varieties} (see \cite{M1-02}, \cite{M2-04}, and \cite{PR-04}).
Moreover, we can use amoebas as an intermediate link between the classical and the tropical geometry. 

The amoeba $\mathscr{A}$  of an algebraic variety $V\subset (\mathbb{C}^*)^n$ is  a  closed subset of $\mathbb{R}^n$, and its (Lebesgue) volume is well-defined.
Passare and Rullg\aa rd  \cite{PR-04}  proved that the area of complex plane curve amoebas is finite and the bound is given in terms of the Newton polygon.
In this paper,
we prove   that the amoeba area  of any  algebraic curve in $(\mathbb{C}^*)^n$ is finite (the area here is with respect  to  the induced Euclidean metric of $\mathbb{R}^n$). 
Moreover, we generalize our result, for any algebraic variety $V$ of dimension $k$ in the  complex algebraic torus $(\mathbb{C}^*)^{2k+m}$  with $m\geq 0$. Our main result is the following theorem: 


\begin{theorem}\label{main theorem}
Let $V$ be an algebraic variety  of dimension $k$  in the complex torus $(\mathbb{C}^*)^{2k+m}$ with $m\geq 0$, such that no  irreducible component of $V$ is contained in a subtorus of dimension less than $2k$. Then, the volume of its amoeba is finite.
\end{theorem} 

Note that if $V$ contains an irreducible  component in an algebraic subtorus of dimension strictly less than $2k$, then the volume of its amoeba is infinite.

\vspace{0.2cm}

The remainder of this paper is organized as follows. In Section 2, we review some properties of the amoebas and also the structure theorem  of the logarithmic limit set defined by Bergman \cite{B-71}, which is used as a tool in the proof of our results. In Section 3,  we prove our main result for complex algebraic curves in $(\mathbb{C}^*)^n$ for any $n\geq 2$. In Section 4, we give an estimate of the bound for algebraic rational complex curves and some examples. In Section 5, we  prove  the main theorem of this paper.

\section{Preliminaries}

Let $V$ be an algebraic variety in $(\mathbb{C}^*)^n$. The {\em amoeba}  $\A$ of $V$ is by definition (see I.M.~Gel'fand, M.M.~Kapranov
 and A.V. Zelevinsky \cite{GKZ-94}) the image of $V$ under the map
\[
\begin{array}{ccccl}
\Log:& &(\mathbb{C}^*)^n&\longrightarrow&\mathbb{R}^n\\
&&(z_1,\ldots ,z_n)&\longmapsto&(\log |z_1|,\ldots ,\log|
z_n|).
\end{array}
\]

Passare and Rullg\aa rd  prove the following (see \cite{PR-04}):

\begin{theorem}\label{Passare-Rullgard}
 Let $f$ be a Laurent polynomial in two variables. Then the area of the amoeba of the curve  with defining polynomial $f$ is not greater than $\pi^2$  times the area of the Newton polytope of $f$.
\end{theorem}

In \cite{MR-00}, Mikhalkin and Rullg\aa rd show that up to multiplication by a constant in $(\mathbb{C}^*)^2$,  the algebraic plane curves with Newton polygon $\Delta$ with  maximal amoeba area are defined over $\R$. Furthermore,  their real loci are isotopic to so-called Harnack curves (possibly singular with ordinary real isolated double points).

\vspace{0.2cm}

Recall that the Hausdorff distance between two closed subsets $A,
B$ of a metric space $(E, d)$ is defined by
$$
d_{\mathcal{H}}(A,B) = \max \{  \sup_{a\in A}d(a,B),\sup_{b\in B}d(A,b)\},
$$
where we take $E =\mathbb{R}^n\times (S^1)^n$, with the distance
defined as the product of the
Euclidean metric on $\mathbb{R}^n$ and the flat metric on $(S^1)^n$.

\vspace{0.2cm}


The {\em logarithmic limit set}  of  a complex algebraic variety  $V$, denoted by $\mathscr{L}^{\infty}(V)$, is the set of limit points of $\mathscr{A}$ in the sphere $S^{n-1} = (\mathbb{R}^n)^*/\mathbb{R}_+$. In other words, if $S^{n-1} $ denotes the boundary of the unit ball $B^n$ and $r$ the map defined by
\[
\begin{array}{ccccl}
r:& &\mathbb{R}^n&\longrightarrow&B^n\\
&&x&\longmapsto&r(x) = \frac{x}{1+|x|} ,
\end{array}
\]
then $ \mathscr{L}^{\infty}(V)= \overline{r(\mathscr{A})}\cap S^{n-1}$. Bergman  \cite{B-71} proved that if $V\subset (\mathbb{C}^*)^n$ is an algebraic variety of dimension $k$, then the cone over  $\mathscr{L}^{\infty}(V)$ is contained in a finite union of $k$-dimensional subspaces of $\mathbb{R}^n$ defined over $\mathbb{Q}$. On the other hand, Bieri and Groves \cite{BG-84} proved that this cone is a finite union of rational polyhedron convex cones of dimension at most $k$, and the maximal dimension  in this union is achieved by at least one polyhedron $P$ in this union. Moreover, one has $\dim_{\mathbb{R}}\mathscr{L}^{\infty}(V) = \dim_{\mathbb{C}}V - 1$. More precisely,  we have  the following structure theorem:

\begin{theorem}
[Bergman, Bieri-Groves] The logarithmic limit set $\mathscr{L}^{\infty}(V)$ of an algebraic variety $V$ in $(\mathbb{C}^*)^n$ is a finite union of rational spherical polyhedrons. The maximal dimension of a polyhedron in this union is achieved at least by one polyhedron $P$ in this union, and we have $\dim_{\mathbb{R}}\mathscr{L}^{\infty}(V) = \dim_{\mathbb{R}}P = \dim_{\mathbb{C}}V - 1$.
\end{theorem}

\section{Area of complex  algebraic curve amoebas}\label{area rational map}

The main result of this section is the following theorem:

\begin{theorem}\label{Curvetheorem}
Let $\mathcal{C}\subset (\mathbb{C}^*)^n$ be an algebraic curve   with $n\geq 2$. Assume that no irreducible component of $\mathcal{C}$ is contained in a subtorus of dimension less than $2$. Then, the area of its amoeba is finite.
\end{theorem}

We start by  proving Theorem \ref{Curvetheorem} in the rational curve case (see Theorem \ref{Area propo}). Recall that  a complex algebraic curve is contained in a subtorus of dimension one, meaning that the curve is the subtorus of dimension one itself (sometimes called a holomorphic annulus). Moreover,  its amoeba is a straight  line in $\mathbb{R}^n$, and this case is not interesting for us because it is not generic (i.e., the Jacobian matrix  of the logarithmic map is not of maximal rank).

Let $n$ and $k$ be two positive integers such that $2k\leq n$.
Let $f: \C^k\longrightarrow (\C^*)^n$ be a rational map, and $V$ be the variety in $(\C^*)^n$ defined by the image of $f$. We denote by $\{z_j\}_{1\leq j\leq k}$ the complex coordinates of $\C^{k}$, and by  $\{f_j\}_{1\leq j\leq n}$ the coordinates of $f$ in $\C^n$. For simplicity, we denote by $\Log f$ the composition $\Log\circ f$. \\

Let $\A_f$ be the amoeba of $V$ (i.e. $\A_f=\Log (V) $). Let $S$ be the set of  points in $\C^{k}$ defined by

\begin{equation*}
S=\{z\in\C^{k}\; \vert\; \textrm{rank}\;\di_z\Log f<2k\}
\end{equation*} 

and $\Si_f=\Log f(S)$ be the set of  critical values of   $\Log f$. \\

By construction, $\Log f$ is an immersion from $\C^{k}\setminus S$ to $\R^n$. Hence, the set $\A_f\setminus\Si_f= \Log f(\C^{k}\setminus S)$ is a  real $2k$-dimensional  immersed submanifold in $\R^n$. We endow $\A_f\setminus \Si_f$ with the induced Riemannian metric $\imath^*\mathcal E_n$, where $\mathcal E_{n}$ is the Euclidean metric of $\R^n$ and $\imath: \A_f\setminus\Si_f\hookrightarrow \R^n$ is the inclusion map. Let $U_f\subset \C^{k}\setminus S$ be an open set such that $\Log f\vert_{U_f}$ is an injective immersion and $\Log f(U_f)=\A_f\setminus\Si_f$. We claim that

\begin{equation}\label{equality volume}
 vol (\A_f\setminus\Si_f, \imath^*\mathcal E_n)=vol(U_f,(\Log f)^*\mathcal E_n),
\end{equation}
where $vol(\A_f\setminus\Si_f, \imath^*\mathcal E_n)$ is the volume of $\A_f\setminus\Si_f$ with respect to the metric $\imath^*\mathcal E_n$. Let $\psi_{2k}$ be a real $2k-$vector field in $\Lambda^{2k}\C^{k}$ which does not vanish on $\C^{k}$ and  $\di v_{(\Log f)^*\mathcal E_n}$,  $\di v_{\mathcal E_{2k}}$ be the volume forms  defined over $U_f$ associated to the metrics $(\Log f)^*\mathcal E_n$ and $\mathcal E_{2k}$ respectively. These two forms are related by the following formula:
\begin{equation}\label{relation volume}
 |\psi_{2k}|_{\mathcal E_{2k}} \di v_{(\Log f)^*\mathcal E_n}=|\psi_{2k}|_{(\Log f)^*\mathcal E_n}\di v_{\mathcal E_{2k}}.
\end{equation}
Now we choose $\psi_{2k}$ such that $\di v_{\mathcal E_{2k}}(\psi_{2k})=|\psi_{2k}|_{\mathcal E_{2k}}=1$. From \eqref{equality volume} and \eqref{relation volume} we deduce that

\begin{equation}\label{volume formula}
 vol(\A_f\setminus\Si_f)=\int_{U_f} \bigl|\psi_{2k}\bigr|_{(\Log f)^*\mathcal E_n}    \di v_{\mathcal E_{2k}},
\end{equation}
and $area:=vol$ if $k=1$. This definition of the volume does not depend on the choice of coordinates. It is more convenient to use the following integral $vol_{2k}$ defined as

\begin{equation}\label{volume multiplicity}
vol_{2k}(\A_f\setminus\Si_f)=\int_{\C^k-S} \bigl|\psi_{2k}\bigr|_{(\Log f)^*\mathcal E_n}    \di v_{\mathcal E_{2k}}.
\end{equation}

\begin{remark}\label{remark 2 vol} $vol_{2k}(\A_f\setminus\Si_f)$ can be viewed as a weighted volume of $\A_f\setminus\Si_f$. In fact, we know that there exist a positive integer $m$ and a family of open connected components on $\A_f\setminus\Si_f$, denoted by $\{R_\alpha\}$, such that

$$\A_f\setminus\Si_f=\bigsqcup_{\alpha=1}^mR_\alpha,\quad   (\Log f)^{-1}R_\alpha=\bigsqcup_{\beta=1}^{p_\alpha}U_{\alpha\beta},\quad\bigsqcup_{\alpha=1}^m\bigsqcup_{\beta=1}^{p_\alpha}U_{\alpha\beta}=\C^{k}\setminus S$$
and the map $\Log f: U_{\alpha\beta}\longrightarrow R_\alpha $ is a diffeomorphism, for all $\alpha\leq m$ and $1\leq \beta\leq p_{\alpha}$.\\

For a chosen  $\beta$, we set $U_f=\bigsqcup_{\alpha=1}^m U_{\alpha\beta}$.
This yields 
\begin{gather}
vol_{2k}(\A_f\setminus\Si_f)=\sum_{\alpha=1}^m p_\alpha\int_{U_{\alpha\beta}} \bigl|\psi_{2k}\bigr|_{(\Log f)^*\mathcal E_n}    \di v_{\mathcal E_{2k}}\label{vol weight},\\
vol(\A_f\setminus\Si_f)=\sum_{\alpha=1}^m \int_{U_{\alpha\beta}} \bigl|\psi_{2k}\bigr|_{(\Log f)^*\mathcal E_n}    \di v_{\mathcal E_{2k}}\label{vol weight 2}.
\end{gather}
If we define $p=\displaystyle\min_{1\leq \alpha\leq m} p_\alpha$ and $P=\displaystyle\max_{1\leq \alpha\leq m} p_\alpha$, then 
\begin{equation}\label{inequality two volumes}
 \frac{vol_{2k}(\A_f\setminus\Si_f)}{P}  \leq vol(\A_f\setminus\Si_f)\leq  \frac{vol_{2k}(\A_f\setminus\Si_f)}{p}.
\end{equation}
\end{remark}

\vspace{0.7cm}

On $\C^{k}$, $\di v_{\mathcal E_{2k}}$ and $\psi_{2k}$ are given by
\begin{equation*}
\di v_{\mathcal E_{2k}}=i^k\di z\wedge\di\bar z,\qquad \psi_{2k}=(-i)^k\frac{\partial}{\partial z}\wedge \frac{\partial}{\partial{\bar z}},
\end{equation*}
where $\di z=\di z_1\wedge\cdots \wedge\di z_k $ and $\frac{\partial}{\partial z}=\frac{\partial}{\partial z_1}\wedge\cdots \wedge\frac{\partial}{\partial z_k}$.

Let us now compute $\bigl|\psi_{2k}\bigr|_{(\Log f)^*\mathcal E_n}$. We have
\begin{equation*}
\frac{\partial \Log f}{\partial z}\wedge \frac{\partial\Log f}{\partial{\bar z}}=\sum_{I=\{i_1<\cdots<i_{2k}\}}\det (\partial_j\Log f_I)_{1\leq j\leq 2k}e_I,
\end{equation*}
where for all $I=\{i_1<\cdots<i_{2k}\}\subset \{1,\ldots,n\}$, $f_I=(f_{i_1},\ldots,f_{i_{2k}})$, $\{e_j\}_{1\leq j\leq n}$ is an orthonormal basis of $T^*\R^n$ and $e_I:=e_{i_1}\wedge\cdots\wedge e_{i_{2k}}$. This implies that 
$$\{e_I\,\vert\,I\subset \{1,\ldots,n\}, \# I=2k\}$$ 
is an orthonormal basis of $\wedge^{2k}T^*\R^n$ with respect to the Euclidean metric. We denote by $\partial_j=\frac{\partial}{\partial z_j}$ if $j\leq k$ and by $\partial_j=\frac{\partial}{\partial \bar z_{j-k}}$ if $j\geq k+1$.   Then

\begin{equation}\label{det g general}
\bigl|\psi_{2k}\bigr|^2_{(\Log f)^*\mathcal E_n}=\sum_{I=\{i_1<\cdots<i_{2k}\}} |\det (\partial_j\Log f_I)_{1\leq j\leq 2k}|^2.
\end{equation}
Hence, we deduce the following inequality:

\begin{equation}\label{ineq vol}
\bigl|\psi_{2k}\bigr|_{(\Log f)^*\mathcal E_n}\leq \sum_{I=\{i_1<\cdots<i_{2k}\}} \big|{\det} (\partial_j\Log f_I)_{1\leq j\leq 2k}\bigr| .
\end{equation}

We have the following result: 
\begin{theorem}\label{Area propo}
Let $f: \C \longrightarrow (\C^*)^n$ be a rational map, and $\mathcal{C}$ be the rational curve of  $(\C^*)^n$ defined by the image of $f$. The area of $\A_f$ with respect to the Euclidean metric of $\R^n$ is finite.
\end{theorem}

\begin{lemma}\label{det g calcul}
Let $f=(f_1,f_2)$ be a rational map from $\C$ to $\C^2$. The function $i\det (\partial_z\Log f, \partial_{\bar z} \Log f )$ is a real-valued rational function. Moreover,
\begin{enumerate}
\item[$(i)$] The map $f$ has simple poles. 
\item[$(ii)$] There exist $P,Q\in\R[X,Y]$ such that $i\det (\partial_z\Log f, \partial_{\bar z} \Log f )=\frac{P}{Q}$ with $\deg\, Q\geq\deg P+3$. 
\end{enumerate}
\end{lemma}

\begin{proof}
We have
$$i\det (\partial_z\Log f, \partial_{\bar z} \Log f )=\frac{i}{4}\biggl(\frac{f_1'}{f_1}\frac{\bar f'_2}{\bar f_2}-\frac{\bar f'_1}{\bar f_1}\frac{f'_2}{f_2}\biggr) .$$

It is trivial that this function, which is $i$ times the Jacobian determinant of $\Log f$, is a real rational function. Its poles are the zeros and the poles of $f$ and their order is  equal to one (even if there is a common pole or zero of $f_1$ and $f_2$, one can check that this pole is also simple). On the other hand, if $i\det (\partial_z\Log f, \partial_{\bar z} \Log f )=\frac{P}{Q}$, then $\de Q\geq\de P+2$. This inequality can be improved. In fact, by elementary computations on the $\de P$ and $\de Q$, we show that $\deg Q\geq\deg P+3$.

\end{proof}

\begin{proof}[Proof of Theorem \ref{Area propo}]

First of all, we do not have to worry about the area of $\Si_f$. Indeed, by Sard's theorem, we know that this area  is equal to  zero. If $\A_f\setminus\Si_f$ is empty, then the area of $\A_f$ is zero. From now on, we assume that $\A_f\setminus\Si_f$ is not empty. Hence, it is a surface defined in $\R^n$ and $\A_f\setminus\Si_f$, $\A_f$ have the same area. 

The area of $\A_f\setminus\Si_f$ is given by \eqref{volume formula}. Hence, it is sufficient to prove that $\bigl|\psi_{2}\bigr|$  is integrable over $\C$.  Inequality \eqref{ineq vol} implies that 

\begin{equation}\label{ineq vol m1}
\bigl|\psi_{2}\bigr|_{(\Log f)^*\mathcal E_n}\leq \sum_{1\leq j<k\leq n} \big|\det (\partial_z\Log f_j, \partial_{\bar z} \Log f_k )\bigr| .
\end{equation}

We claim that all the functions on the right-hand side of \eqref{ineq vol m1} are integrable. Indeed, using Lemma \ref{det g calcul}, $(i)$, we have the integrability in a neighborhood of any pole.  By $(ii)$, we get the integrability at infinity.  
\end{proof}

\vspace{0.3cm}

\begin{proof}[Proof of Theorem \ref{Curvetheorem}]
Let $\mathcal{C}$ be an algebraic curve in $(\mathbb{C}^*)^n$. Then its closure $\overline{\mathcal{C}}$ in $\mathbb{C}\mathbb{P}^n$ is an algebraic curve. Hence,  any end of $\mathcal{C} $  corresponds to a local branch of $\overline{\mathcal{C}}$  at  some point $p\in\partial \overline{\mathcal{C}}\setminus \mathcal{C}$. After a monomial map of $(\mathbb{C}^*)^n$ if necessary, we may assume that  $p$ corresponding to an end of $\mathcal{C} $ is the origin of $\mathbb{C}^n$ in $\mathbb{C}\mathbb{P}^n$.

A local parametrization $\rho_p$ of  a branch of  $\overline{\mathcal{C}}$ at $p$ can be written in terms of vectorial  Puiseux series  in $t$ near zero as follows:
\[
\begin{array}{cccl}
\rho_p: &\mathbb{C}^*&\longrightarrow & (\mathbb{C}^*)^n\\
& t &\longmapsto & (b_1t^{u_1},\ldots ,b_nt^{u_n}),
\end{array}
\]
where $(b_1,\ldots ,b_n)\in (\mathbb{C}^*)^n$, and $(u_1,\ldots ,u_n)\in \mathbb{Q}_{\geq 0}^n$.
Indeed,  the  Bergman logarithmic limit set of a curve is a finite number of points $\{v_i\}$ on the sphere $S^{n-1}$ (see \cite{B-71}). By Bieri and Groves (see \cite{BG-84}),
 if $O$ denotes the origin of  $\mathbb{R}^n$, then the slope $\vec{u}$ of the real line $(Ov_i)$ in $\mathbb{R}^n$ is rational (the slope here means the direction vector of the line, and rational means that its coordinates are rational). Hence, there exist real lines $L_{\vec{u},\, j}$ in $\mathbb{R}^n$ parametrized by
 $$
  x \longmapsto (a^{j}_1+ xu_1,\ldots ,a^{j}_n+ xu_n)
 $$
 with $x\in \mathbb{R}$ and $(a^{j}_1,\ldots ,a^{j}_n)\in \mathbb{R}^n$, such that the amoebas $\mathscr{A}_h$  of $h$ in the defining ideal of the curve $\mathcal{C}$ reaches all these lines at the infinity in the direction $\vec{u}$. So, for each line $L_{\vec{u},\, j}$ there exists $(b_1,\ldots ,b_n)\in(\mathbb{C}^*)^n$ such that the Hausdorff distance between the complex line $\mathcal{L}_{\vec{u},\, j}$ parametrized by
 $$
  t \longmapsto (b_1t^{u_1},\ldots ,b_nt^{u_n}),
 $$
 and $V_h\cap \Log^{-1}(L_{\vec{u},\, j}\setminus B(O, R))^+$ tends to zero when $R$ goes to infinity, where $(L_{\vec{u},\, j}\setminus B(O, R))^+$ denotes the component which is in the direction of $\vec{u}$. Hence, the Hausdorff distance between the intersection of the curve with  $\Log^{-1}(L_{\vec{u},\, j}\setminus B(O, R))^+$ and $\mathcal{L}_{\vec{u},\, j}$ tends to zero for $R$ sufficiently large. Now, using Theorem \ref{Area propo} and the fact that the number of ends of an algebraic curve is finite, we obtain the result.
\end{proof}

\section{An estimate for the area of rational curve amoebas}

In this section, we assume that $k=1$. Recall that $f$ is  the rational map defined in Section \ref{area rational map} with $k=1$. Hence, for any integer  $j\in[1,n]$:
$$f_j(z)=c\prod_{\ell=1}^{d_j} (z-a_{j\ell})^{m_{j\ell}},$$ 
where $a_{j\ell}$ are distinct poles or zeros of $f_j$, $m_{j\ell}$ are their multiplicities ($m_{j\ell}$ are negative in the case of poles) and $d_j$ is the number of distinct zeros and poles of $f_j$. We define the positive integers $n_j:=\sum_{\ell=1}^{d_j}|m_{j\ell}|$. These integers represent the number of poles and zeros of $f_j$ counted with their multiplicities.
\begin{theorem}\label{proposition inequalities}
Let $p$ be the positive integer defined by 
$$p= \min_{x\in \A_f\setminus\Si_f}\#(\Log f)^{-1}x.$$ 
The following inequalities always hold:  
\begin{gather}\label{inequality area fond}
p \cdot area (\A_f)\leq vol_{2} (\A_f)\leq \pi^2\sum_{1\leq j_1<j_2\leq n}n_{j_1}n_{j_2}.
\end{gather}
Moreover, $p \cdot area (\A_f)=vol_{2} (\A_f)$ if and only if $\Log f: \C^k\setminus S\longrightarrow \A_f\setminus\Si_f$ is a covering with exactly p sheets.
\end{theorem}

\begin{proof}
Recall that $vol_{2}$ is defined by \eqref{volume multiplicity}.\\
If ${n= 2}$, using \eqref{volume multiplicity}, \eqref{det g general}, we obtain
\begin{equation*}
vol_{2}(\A_f)\leq \sum_{\ell=1,\ell'=1}^{d_1,d_2}|m_{1\ell}m_{2\ell'}|vol_{2}(\A_{f_{a_{1\ell},a_{2\ell'}}}),
\end{equation*}
with $f_{a_{1\ell},a_{2\ell'}}(z)=(z-a_{1\ell},z-a_{2\ell})$. We know that $vol_{2}(\A_{f_{a_{1\ell},a_{2\ell'}}})=\pi^2$ (it can be proven, using the substitution $z=(a_{1\ell}-a_{2\ell})t+a_{1\ell}$ and Example 1 below for $m=1$, which is a plane line). Hence,
\begin{equation}\label{inequality n2}
vol_{2}(\A_f)\leq\pi^2\sum_{\ell=1,\ell'=1}^{d_1,d_2}|m_{1\ell}m_{2\ell'}|.
\end{equation}
If $n\geq 2$, using \eqref{ineq vol} and integrating, we obtain
\begin{equation}\label{cccf}
vol_{2} (\A_f)\leq \sum_{1\leq j_1<j_2\leq n}vol_{2}(\A_{(f_{j_1},f_{j_2})}).
\end{equation}
A combination of \eqref{inequality n2} and \eqref{cccf} yields

\begin{equation}
vol_{2} (\A_f)\leq \pi^2\sum_{1\leq j_1<j_2\leq n}\sum_{\ell=1,\ell'=1}^{d_{j_1},d_{j_2}}|m_{j_1\ell}m_{j_2\ell'}|.
\end{equation}
Hence
\begin{equation}
vol_{2}(\A_f)\leq \pi^2\sum_{1\leq j_1<j_2\leq n}n_{j_1}n_{j_2},
\end{equation}
which gives the second inequality of Theorem \ref{proposition inequalities}. The first one is a consequence of \eqref{inequality two volumes}. \\
Assume that $p \cdot area (\A_f)=vol_{2} (\A_f)$. Using \eqref{vol weight}, \eqref{vol weight 2}, we obtain that $p_\alpha=p$, for all $\alpha\leq m$. This means that the number of connected components in $(\Log f)^{-1} R_\alpha$ does not depend on $\alpha$ and $\Log f$ is a $p-$sheeted covering map.\\
If we suppose that $\Log f: \C^k\setminus S\longrightarrow \A_f\setminus\Si_f$ is a covering with exactly $p$ sheets, then \eqref{inequality two volumes} becomes an equality. 
\end{proof}

\subsection*{Examples}

\begin{enumerate}
\item\label{example 1} Let us compute the area of the amoeba $\A_f$ with $f(z)=(z,z^m-1)$ and $m\in\mathbb Z^*$. The set of singular points $S$ is the union of the half lines given by
$$S=\bigcup_{j=1}^{2|m|}\biggl\{z\in\C\,\vert\,\arg z=\frac{j\pi}{|m|}\biggr\}.$$
The set of critical values is $\Si_f=\Log f(S)$, which is a curve in $\R^2$ and bounds $\A_f$. The map $\Log f$ is a covering map with exactly $2|m|$ sheets and for any integer $j\in[1,2|m|]$ we have that
$$\Log f:U_j=\biggl\{z\in\C^*\,\vert\,\frac{(j-1)\pi}{|m|}<\arg z<\frac{j\pi}{|m|}\biggr\}\longrightarrow \A_f - \Si_f$$
is a diffeomorphism. Hence, $\Log f$ is a $2|m|-$sheeted covering map. By \eqref{volume formula} and \eqref{det g general} we obtain 

\begin{equation*}
area(\A_{f})=\int_{U_1}\biggl|\det (\frac{\partial\Log f}{\partial z},\frac{\partial \Log f}{\partial \bar z})\biggr|i\di z\wedge\di\bar z=\int_{U_1}\frac{|m||z^m-\bar z^m|}{4|z|^2|z^m-1|^2}i\di z\wedge\di\bar z=\frac{\pi^2}{2|m|} ,
\end{equation*}
with $area =vol$. We deduce that $vol_2(\A_f)=2|m| area(\A_f)=\pi^2$.
\item Now, let us consider the real line in $(\mathbb{C}^*)^3$ parametrized by $g(z)=(z,\\
z+\frac{1}{2},z-\frac{3}{2})$. The amoeba $\A_g$ is a surface  in $\R^3$ with boundary as we can see in  Figure 1 (this fact is proven in \cite{NP-10}).  The set of singular points is the line of real points,  and $\Log g$ is  a 2-sheeted covering map. It is complicated to compute the area of $\A_g$. However, using the estimate given in Theorem \ref{proposition inequalities}, we deduce that $2\,area (\A_g)=vol_2(\A_g)\leq 3\pi^2$.
\item Let $h(z)=(z,z+1,z-2i)$ be the parametrization of a complex line in $(\C^*)^3$. The amoeba $\A_h$ is a surface  in $\R^3$ without boundary as we can see in Figure 2, and topologically it is a Riemann sphere with four marked points.
Note that this line is not real and the set of critical values of $\Log$ restricted to this line is empty. The map $\Log h:\; \C-\{-1,0,2i\}\longrightarrow \A_h$ is a diffeomorphism. By Theorem \ref{proposition inequalities}, $area (\A_h)=vol_2(\A_h)\leq 3\pi^2$.
\end{enumerate}

\begin{figure}[h!]
\begin{center}
\includegraphics[width=1\textwidth]{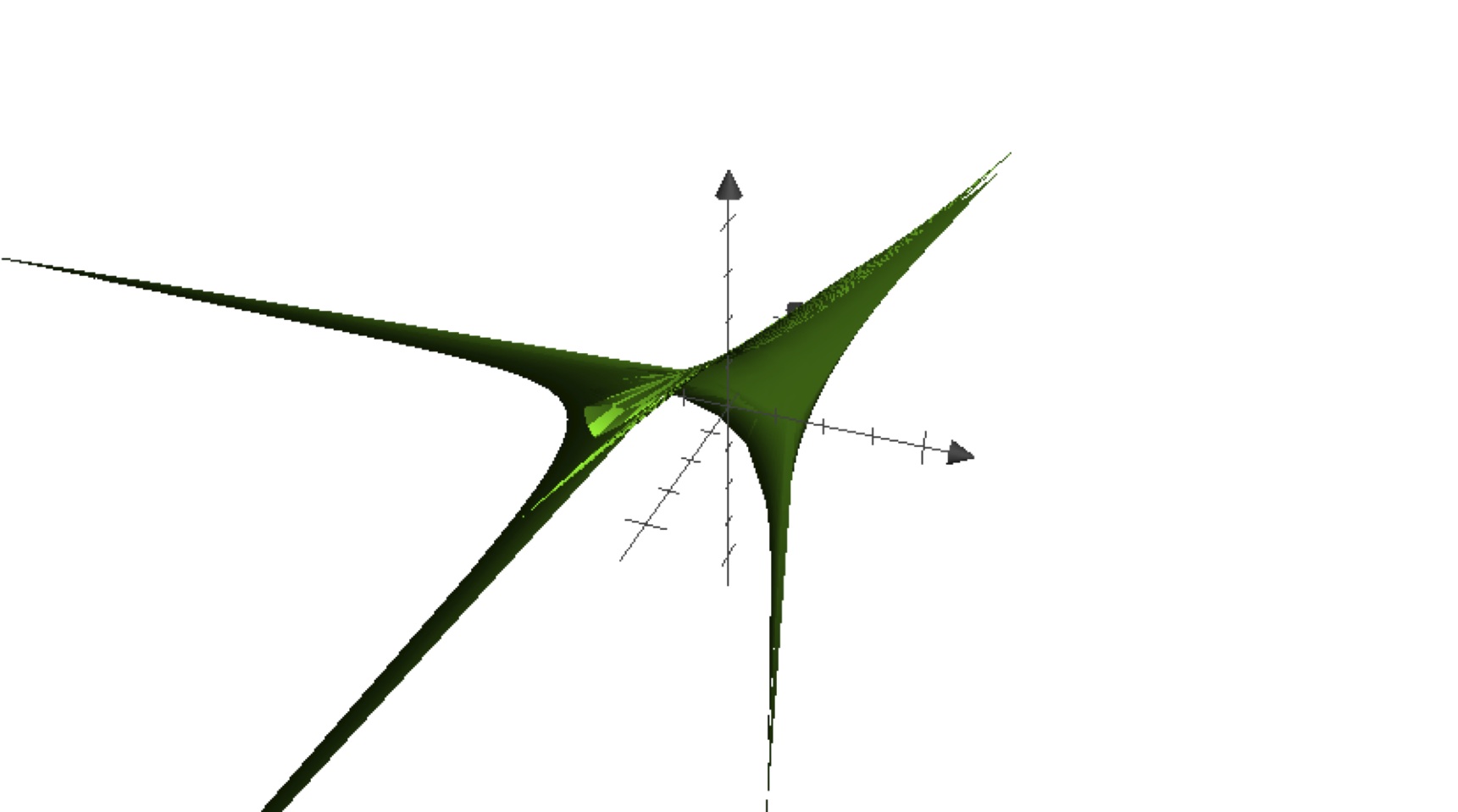}
\caption{The amoeba of the real line in $(\mathbb{C}^*)^3$ given by the parametrization $g(z)=(z,z+\frac{1}{2},z-\frac{3}{2})$.}
\label{c}
\end{center}
\end{figure}

\begin{figure}[h!]
\begin{center}
\includegraphics[width=1\textwidth]{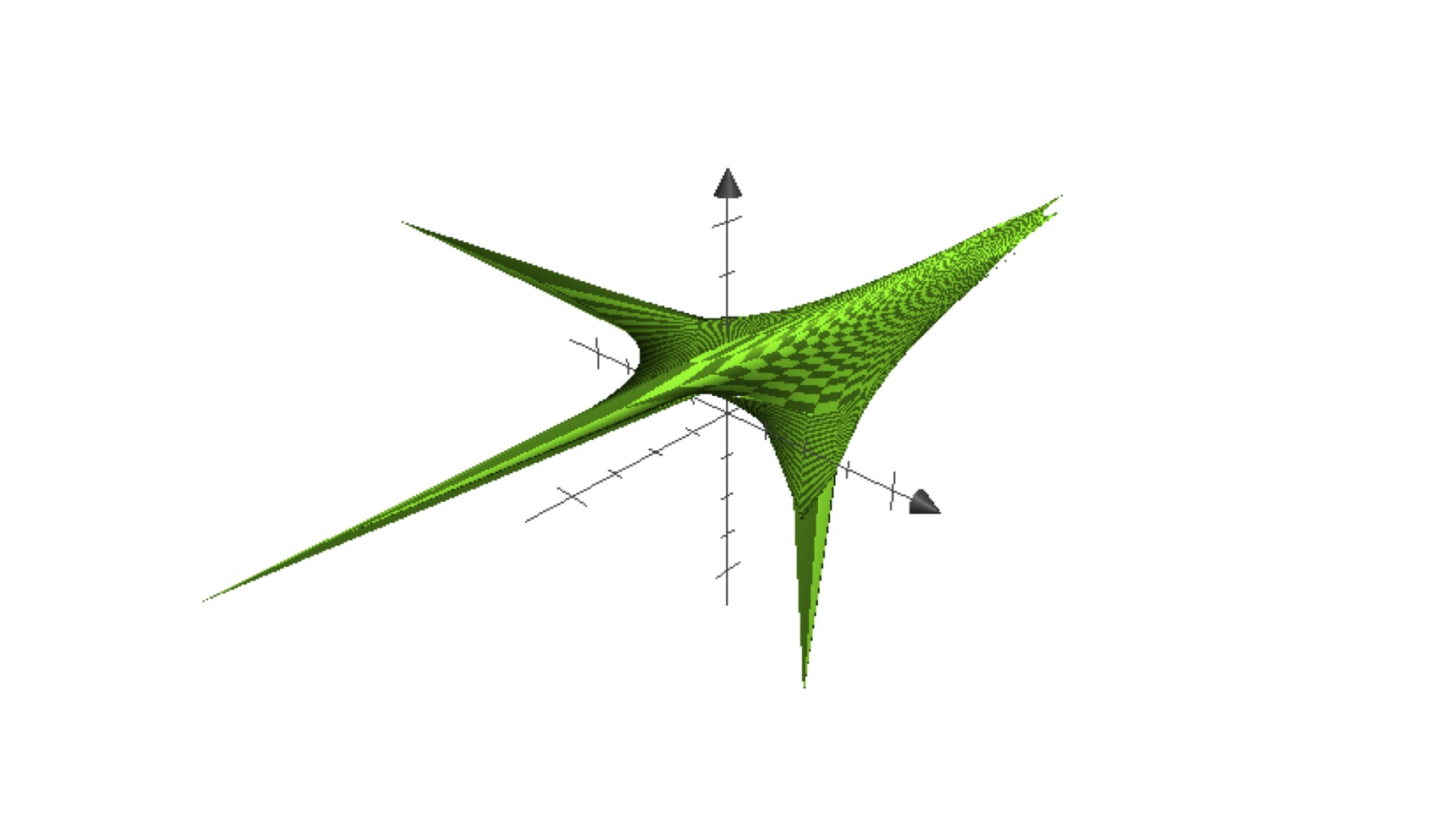}
\caption{The complex line in $(\mathbb{C}^*)^3$ given by the parametrization $h(z)=(z,z+1,z-2i)$.}
\label{c}
\end{center}
\end{figure}

\vspace{0.3cm}

Passare and Rullg\aa rd give an estimate for amoeba areas of  complex algebraic  plane curves (see Theorem \ref{Passare-Rullgard}). Our estimate works only for rational curves immersed in $(\C^*)^n$. However, in this case we have  a finer estimate (see Theorem~\ref{proposition inequalities}). Indeed, if we consider Example \ref{example 1}, $area (\A_{f})= \frac{\pi^2}{2m}$, $n_1=1$, $n_2=m$, $p=2m$ and the area of Newton's polygon is $\frac{m}{2}$. Inequality \eqref{inequality area fond} gives $area (\A_{f})\leq \frac{\pi^2}{2}$, and the Passare-Rullg\aa rd estimate gives    $area (\A_{f})\leq \frac{m\pi^2}{2}$.

\vspace{0.2cm}

\section {Volume of  generic complex algebraic  variety amoebas}

In this section, we assume that $n=2k+m$ is an integer with $k\geq 1$ and $m\geq 0$. Let $V\subset (\mathbb{C}^*)^n$ be an algebraic variety of dimension $k$ with defining ideal $\mathcal{I}(V)$ and $\mathscr{L}^{\infty}(V)$ be its logarithmic limit set. We denote by $\Verte (\mathscr{L}^{\infty}(V))$ the set of vertices of $\mathscr{L}^{\infty}(V)$. Let $V\subset (\mathbb{C}^*)^n$ be a generic  algebraic variety of dimension $k$. Let $v\in \Verte (\mathscr{L}^{\infty}(V))$ and let us denote by 
$\mathscr{H}_{R,\, v}$ the hyperplane in $\mathbb{R}^n$ with normal the vector $\vec{Ov}$ such that $d(O, \mathscr{H}_{R,\, v}) = R$, where $O$ is the origin of $\mathbb{R}^n$  and  $R\in\mathbb{R}_+$  is sufficiently large. We denote by $\mathscr{H}^-_{R,\, v}$ the half space with boundary $\mathscr{H}_{R,\, v}$  containing the origin.

In this section, we prove Theorem \ref{main theorem}, using the following proposition:
 
 \vspace{0.2cm}

\begin{proposition} \label{prop prop}
 With the above notation, $V\setminus \Log^{-1} (\mathscr{H}^-_{R,\, v})$  is a fibration over an algebraic variety  $V_v$ of dimension $k-1$ contained in $(\mathbb{C}^*)^{n-1}$, and its fibers are the ends of algebraic curves in $(\mathbb{C}^*)^n$. Moreover, these ends have a  rational parametrization with the same slope (i.e., their image under $\Log$ are lines with the same slope).
\end{proposition}

We start by proving the following lemma:
\begin{lemma}
Let $V \subset (\mathbb{C}^*)^n$ be an algebraic  variety of dimension $k$. 
 Then, for each vertex $v$ of its logarithmic limit set $\mathscr{L}^{\infty}(V)$, we have the following: there exists a  complex  algebraic variety $V_v\subset (\mathbb{C}^*)^{n-1}$  of dimension $k-1$ such that the boundary of the Zariski closure  $\overline{V}$  of $V$ in  $(\mathbb{C}^*)^{n-1}\times \mathbb{C}$  is equal to $V_v$ (i.e., $\partial \overline{V} = \overline{V}\setminus V = V_v$, where $(\mathbb{C}^*)^{n-1} = (\mathbb{C}^*)^{n-1}\times \{ 0\}\subset (\mathbb{C}^*)^{n-1}\times \mathbb{C}$).
\end{lemma}

\begin{proof} 
If  $v$ belongs to $\Verte (\mathscr{L}^{\infty}(V))$,  then after a monomial map defined by a matrix $A_v\in GL_n(\mathbb{Z})$ if necessary, we may assume that $v=(0,\ldots ,0,-1)\in \mathbb{S}^{n-1}$.
 Let $\mathbb{C}[z_1^{\pm 1},\ldots ,z_{n-1}^{\pm 1}, z_n]\subset \mathbb{C}[z_1^{\pm 1},\ldots ,z_{n}^{\pm 1}]$ be the inclusion of rings and $ \phi : \mathbb{C}[z_1^{\pm 1},\ldots ,$ \\ $z_{n-1}^{\pm 1},z_n] \rightarrow   \mathbb{C}[z_1^{\pm 1},\ldots ,z_{n-1}^{\pm 1}]$ be the homomorphism which sends $z_n$ to zero. Let $\mathcal{J} = \mathcal{I}(V)\cap \mathbb{C}[z_1^{\pm 1},\ldots ,z_{n-1}^{\pm 1}, z_n]$ and $\mathcal{I}_v$ be its image in $\mathbb{C}[z_1^{\pm 1},\ldots ,z_{n-1}^{\pm 1}]$. We denote by $V_{(n-1)}$ the subvariety of $(\mathbb{C}^*)^{n-1}\times \mathbb{C}$ defined by $\mathcal{J}$, and by $V_v = V_{(n-1)}\cap (\mathbb{C}^*)^{n-1}\times \{ 0\}$ the subvariety defined by $\mathcal{I}_v$. We check that $V_{(n-1)} = \overline{V}$ where $\overline{V}$ denotes the Zariski closure of $V$ in $(\mathbb{C}^*)^{n-1}\times \mathbb{C}$ and $V_v $ is the boundary of $\overline{V}$ i.e., $\partial \overline{V} = \overline{V}\setminus V = V_v $.
\end{proof}

\vspace{0.2cm}

\begin{proof}[Proof of Proposition \ref{prop prop}] For each point $x$ in $V_v$ there exists an algebraic curve $\mathcal{C}_{x}$ in $V$ such that its closure in $(\mathbb{C}^*)^{n-1}\times \mathbb{C}$ contains $x$ and its logarithmic limit set contains the point  $v$. Indeed, we have the following commutative diagram:
\begin{equation}
\xymatrix{
(\mathbb{C}^*)^{n-1}\times\mathbb{C}\ar[d]_{\Log_{|(\mathbb{C}^*)^n}}\ar[rr]^{\pi_{(n-1)}^{\mathbb{C}}}&&(\mathbb{C}^*)^{n-1}\ar[d]^{\Log}\cr
\mathbb{R}^n\ar[rr]^{\pi_{(n-1)}^{\mathbb{R}}}&&\mathbb{R}^{n-1},
}\nonumber
\end{equation}
where $\pi_{(n-1)}^{\mathbb{C}}$ and $\pi_{(n-1)}^{\mathbb{R}}$ are the projections on $(\mathbb{C}^*)^{n-1}$ and $\mathbb{R}^{n-1}$ respectively.
The limit of $\pi_{(n-1)}^{\mathbb{C}}(V\setminus\mathscr{H}_{v,\, R}^-)$ when $R$ goes to infinity is equal to $V_v$ (with respect to the Hausdorff metric on compact sets). Furthermore, the limit of $\pi_{(n-1)}^{\mathbb{R}}(\Log (V\setminus\mathscr{H}_{v,\, R}^-))$ when $R$ goes to infinity is equal to the amoeba of $V_v$. Hence, the limit of $\pi_{(n-1)}^{\mathbb{C}}(\mathcal{C}_{x}\setminus\mathscr{H}_{v,\, R}^-)$ when $R$ goes to infinity contains the point $x$.  The end of $\mathcal{C}_{x}$ corresponding to $v$ and containing $x$ is parametrized as
\[
\begin{array}{cccl}
\rho_v:&\mathbb{C}^*&\longrightarrow&(\mathbb{C}^*)^n\\
& t &\longmapsto&(b_{x,\, 1}t^{u_1},\ldots ,b_{x,\, n}t^{u_n}),
\end{array}
\]
where the coefficients $b_{x,\, j}$ depend on the holomorphic branch of $\overline{\mathcal{C}_x}$ at $x$, and the powers $u_j$ depend only on $v$. Moreover, 
for any $x_1\ne x_2$ in $V_v$ the end of the curve $\mathcal{C}_{x_1}$ corresponding to $x_1$, and the end of the curve $\mathcal{C}_{x_2}$ corresponding to $x_2$ have an empty intersection.  It may be that the curves $\mathcal{C}_{x_1}$ and $\mathcal{C}_{x_2}$ are the same. In fact,  in this case this means that the  curve has more than one end corresponding to $v$.  In other words, if $\overline{\mathcal{C}}_{x_1}$ is the Zariski closure of $\mathcal{C}_{x_1}$ in $\mathbb{CP}^n$, then it 
has more than one holomorphic branch at $v$. Indeed, if the intersection of these ends is not empty, then from the fact that they are holomorphic and with the same slope, they should be equal. This is a  contradiction with  the assumption on $x_1$ and $x_2$. Hence, for $R$   sufficiently large, $V\setminus\mathscr{H}_{v,\, R}^-$ is a fibration over $V_v$.
\end{proof}

\vspace{0.2cm}

\begin{proof}[Proof of Theorem \ref{main theorem}] 
Recall that the volume is always computed  with respect to the induced measure of the ambient  space. 
Using  induction on the dimension $k$ of the variety, Proposition  5.1  and Theorem 4.1,  there exists a rational number $q_v$  depending only on  $v$ and on the variety $V_v$ such that  the inequality
$$
vol (\Log (V\setminus\mathscr{H}_{v,\, R}^-)) \leq \pi^2 q_v vol (\mathscr{A}(V_v))
$$
holds.  There is a finite number of vertices of the logarithmic limit set, so there exists a positive real number $K$ such that the following inequality holds:
$$vol (\mathscr{A}(V)) \leq K + \sum_{v\in\Verte (\mathscr{L}^{\infty}(V))}\pi^2q_v vol (\mathscr{A}(V_v)).$$
\end{proof}

\vspace{0.3cm}

%
%
%
%
%
%


%
%
%
\end{document}